\documentclass{siamart0216}
\usepackage{color,amsmath,amssymb}
\usepackage{graphicx, subfig}
\usepackage[procnames]{listings}

\title{Decomposition of stencil update formula into atomic stages}
\author{Qiqi Wang}

\definecolor{keywords}{rgb}{0.6,0.4,0.0}
\definecolor{comments}{rgb}{0,0,0.4}
\definecolor{red}{rgb}{0.5,0,0}
\definecolor{green}{rgb}{0,0.4,0}
\definecolor{lightgray}{rgb}{0.9,0.9,0.9}
\begin{document}
\maketitle

\begin{abstract}
    In parallel solution of partial differential equations, a complex stencil update formula that accesses multiple layers of neighboring grid points sometimes must be decomposed into atomic stages, ones that access only immediately neighboring grid points.  This paper shows that this requirement can be formulated as constraints of an optimization problem, which is equivalent to the dual of a minimum-cost network flow problem.  An optimized decomposition of a single stencil on one set of grid points can thereby be computed efficiently.
\end{abstract}



\section{Introduction}


In large-scale solution of partial differential equations,
there is a need to better exploit the
computational power of massively parallel systems~\cite{dongarra2014applied,
pi2013cfd, larsson2014prospect}.
Particularly needed are
algorithms that hide network and memory latency~\cite{
chen1992reducing,wonnacott2000using, ghysels2013hiding, strumpen1994exploiting,
datta2009optimization},
increase computation-to-communication ratios~\cite{rostrup2010parallel,
crovella1992using, barney2010introduction, boyd1995modeling}, and
minimize synchronization~\cite{demmel2008avoiding, hoemmen2010communication,
carson2013avoiding}.
It is challenging to design and deploy solutions to address these problems,
partly due to the variety and complexity of schemes commonly used for
solving partial differential equations.  A solution may have been
demonstrated on a simple scheme 
with a compact stencil, but to support complex discretization scheme
involving arbitrary stencil can sometimes be challenging or cumbersome.
Because of this, it would be useful to decompose complex discretization schemes
into components that are easier to integrate into stencil computation
algorithms.

This paper describes a process that decomposes a general update formula
into an optimal series of atomic update formulas, each with a
compact stencil that involves only the immediate neighbors in a mesh.
This process enables application of algorithms and software
that operate only on atomic update formulas to formulas with larger stencils.
The algorithm
can also serve as a pre-processing step for stencil compilers such as
OP2~\cite{mudalige2012op2,reguly2016acceleration},
Patus~\cite{christen2011patus}, and Simit~\cite{kjolstad2016simit},
as well as loop transforming
techniques such as tiling~\cite{renganarayanan2007parameterized}.

\subsection{Stencil update formula and atomic decomposition}
When performing a variety of physical simulations, we discretize space
into grid points, and time into time steps.  The resulting discretized
equation often updates a few values in each grid point every time step,
following a predetermined stencil update formula.
For example, conduction of heat in one-dimensional, homogeneous
structures is often modeled by the 1D heat equation
\begin{equation}
\frac{\partial u}{\partial t} = \frac{\partial^2 u}{\partial x^2}.
\end{equation}
It can be simulated with the update formula
\begin{equation} \label{update00}
u_i^{n+1} = u_i^n + \Delta t \frac{u_{i-1}^n - 2 u_i^n + u_{i+1}^n}{\Delta x^2},
\end{equation}
where subscript $i$ denotes spatial grid point, and superscript $n$ denotes time step.  The set of neighboring grid points involved, $\{i-1, i, i+1\}$, is called the stencil of this update formula.
This update formula is derived through manipulation of Taylor series,
formally by approximating the spatial derivative with a linear combination
of neighboring values, a technique known as finite difference,
in conjunction with a time advancing method called forward Euler.

To automate the manipulation of these update formulas, it is useful to
represent them in a computer language.
In Python, this update formula can be described as the following function
\begin{lstlisting}
def heat(u0):
    return u0 + Dt/Dx**2 * (im(u0) - 2*u0 + ip(u0))
\end{lstlisting}
where \lstinline{im} and \lstinline{ip} represent values at the $i-1$st grid
point and $i+1$st grid point, respectively.
The same stencil update formula (\ref{update00}) is applied at every
grid point $i$, for every time step $n$.

We can solve a wide variety of problems by applying stencil update formulas
at a set of spatial grid points over a series of time steps.
More complex update formulas are often used to increase the accuracy,
and to solve more complex equations.
To increase the accuracy for solving the same 1D heat equation, for example, one
may upgrade the time advancing method from forward Euler to the midpoint
method.  Also known as the second-order Runge-Kutta, it is derived through more
complex manipulation of Taylor series.  The resulting update formula is
\begin{equation} \label{update01}
\begin{split}
u_i^{n+\frac12} &= u_i^n + \frac{\Delta t}2 \frac{u_{i-1}^n - 2 u_i^n + u_{i+1}^n}{\Delta x^2}, \\
u_i^{n+1} &= u_i^n + \Delta t \frac{u_{i-1}^{n+\frac12} - 2 u_i^{n+\frac12} + u_{i+1}^{n+\frac12}}{\Delta x^2}. \\
\end{split}
\end{equation}
Because $u_i^{n+1}$ depends on $u_{i-2}^n$ and $u_{i+2}^n$, the stencil of this update formula is $\{i-2,i-1,i,i+1,i+2\}$.
This update formula can be described as the following function
\footnote{Note that running this descriptive function does not
necessarily perform the computation.  It can merely build a data structure
containing information about the steps required to perform
the computation.}
\begin{lstlisting}
def heatMidpoint(u0):
    u_half = u0 + Dt/Dx**2/2 * (im(u0) - 2*u0 + ip(u0))
    return u0 + Dt/Dx**2 * (im(u_half) - 2*u_half + ip(u_half))
\end{lstlisting}

Scientists and mathematicians invented numerous updating formulas to
simulate various problems.  To make them accurate, stable, flexible,
and appealing in other aspects, they craft formulas that can involve
orders of magnitude more calculations than those in our examples.
In this paper, we focus on update formulas that use a fixed number of
inputs and produce a fixed number of outputs at every grid points.
The outputs depend on the inputs at a stencil, a neighboring set
of grid points.  The update formula is applied to a set of
grid points over a series of time steps.

We can decompose a complex update formula into a sequence of stages.
Such decomposition is desirable if the simulation runs on
massively parallel computers.  Processors in such computers must
communicate during a simulation.  These communications can be simplified if
a complex update formula is decomposed into simpler stages in the following
way
\begin{enumerate}
\item Each stage generates outputs that feed into the inputs of the next stage.
      The inputs of the first stage and the outputs of the last stage match
      the inputs and outputs of the entire update formula.
\item Each stage is {\bf atomic}, which means that its outputs at
    each grid point depend on the
    inputs at no further than the immediately neighboring grid points.
\end{enumerate}
Some parallel computing method, such as the swept decomposition scheme,
is based on the assumption that an update formula is decomposed into atomic
stages.

For example, the update formula (\ref{update00}) is an atomic stage.
Its input is $u_i^n$; its output is $u_i^{n+1}$.  Update formula
(\ref{update01}) can be decomposed into two atomic stages in the following
way.  The input of the first stage is $u_i^n$, and the outputs include 
$u_i^{n+\frac12}$ and a copy of $u_i^{n}$.  These outputs must be the inputs of
the next stage, whose output is $u_i^{n+1}$.  These stages can be encoded as
\begin{lstlisting}
def heatMidpoint_stage1(u0):
    u_half = u0 + Dt/Dx**2/2 * (im(u0) - 2*u0 + ip(u0))
    return u0, u_half

def heatMidpoint_stage2(inputs):
    u0, u_half = inputs
    return u0 + Dt/Dx**2 * (im(u_half) - 2*u_half + ip(u_half))
\end{lstlisting}

Decomposition into atomic stages is not unique.  In addition to
the decomposition above, for example, the same update formula (\ref{update01})
can be decomposed in the following different ways:
\begin{enumerate}
        \item
\begin{lstlisting}
def heatMidpoint_stage1(u0):
    im_plus_ip_u0 = im(u0) + ip(u0)
    return u0, im_plus_ip_u0

def heatMidpoint_stage2(inputs):
    u0, im_plus_ip_u0 = inputs
    u_half = u0 + Dt/Dx**2/2 * (im_plus_ip_u0 - 2*u0)
    return u0 + Dt/Dx**2 * (im(u_half) - 2*u_half + ip(u_half))
\end{lstlisting}

        \item
\begin{lstlisting}
def heatMidpoint_stage1(u0):
    u_half = u0 + Dt/Dx**2/2 * (im(u0) - 2*u0 + ip(u0))
    im_plus_ip_u_half = im(u_half) + ip(u_half)
    return u0, u_half, im_plus_ip_u_half

def heatMidpoint_stage2(inputs):
    u0, u_half, im_plus_ip_u_half = inputs
    return u0 + Dt/Dx**2 * (im_plus_ip_u_half - 2*u_half)
\end{lstlisting}
\end{enumerate}
Note that the first alternative decomposition passes two variables from the
first stage to the second stage, the same number as passed by the original
decomposition; the second alternative decomposition, however, passes three
variables from the first stage to the second stage.  Because passing
variables between stages may incur communication of data between parallel
computing units, we consider it less efficient to pass more variables.
In this metric, the second alternative decomposition is inferior to both
the original decomposition and the first alternative.

Many update formulas can be decomposed into a sequence of atomic stages,
such that the outputs of each stage is the inputs of the next.
The goal of this paper is to automatically find the best decomposition
for very complex update formulas, such that the total amount of variables
passed between the decomposed stages is as few as possible.  After
decomposition, a stencil optimization framework such as Modesto~\cite{gysi2015modesto}
can further improve the performance of the stages.

\subsection{Motivation by the swept rule}
The algorithm developed in this paper
is motivated mainly the the swept
rule~\cite{alhubail2016swept, alhubail2016swept2d}.
It is an algorithm
designed to break the latency barrier of parallel computing by minimizing
synchronization.  This algorithm has the potential to significantly
increase the strong scaling limit~\cite{fischer2015scaling} of many simulations,
allowing them to better exploit the computational power of massively parallel
systems.

The algorithm has been demonstrated on atomic
stencil update formulas, i.e., ones with compact stencils, involving only
the immediate neighbor in a mesh.  Although it was suggested that a more
complex formula can be decomposed into a series of atomic update formulas,
it is not obvious how to effectively do so in general.
An algorithmic approach to performing such decomposition is desirable.

\section{Graph theoretical representation of a stencil update formula and its
atomic decomposition}
\label{sec_graph}

\begin{figure}[htb!] \centering
    \begin{lstlisting}
    def heatMidpoint(u0):
        uHalf = u0 + Dt/Dx**2/2 * (im(u0) - 2*u0 + ip(u0))
        return u0 + Dt/Dx**2 * (im(uHalf) - 2*uHalf + ip(uHalf))
    \end{lstlisting}
    \includegraphics[width=\textwidth]{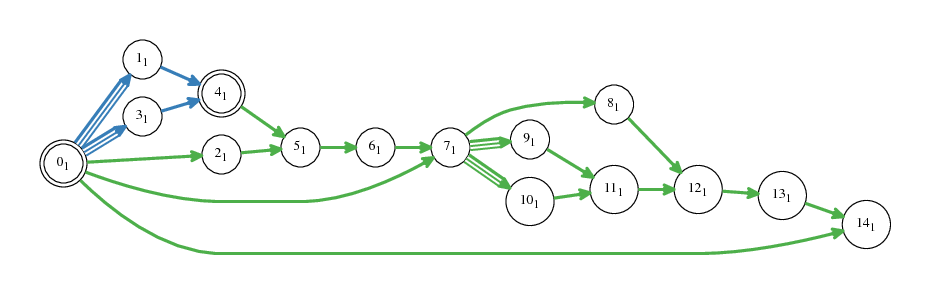} \\
    \begin{tabular}{|c|c|}
        \hline
        Vertex & represents intermediate value \\
        \hline
        0 & \lstinline!u0! \\
        1 & \lstinline!im(u0)! \\
        2 & \lstinline!2 * u0! \\
        3 & \lstinline!ip(u0)! \\
        4 & \lstinline!im(u0) + ip(u0)! \\
        5 & \lstinline!im(u0) + ip(u0) - 2 * u0! \\
        6 & \lstinline!dt / dx**2 * (im(u0) + ip(u0) - 2 * u0)! \\
        7 & \lstinline!uHalf = u0 + dt / dx**2 * (im(u0) + ip(u0) - 2 * u0)! \\
        8 & \lstinline!2 * uHalf! \\
        9 & \lstinline!im(uHalf)! \\
       10 & \lstinline!ip(uHalf)! \\
       11 & \lstinline!im(uHalf) + ip(uHalf)! \\
       12 & \lstinline!im(uHalf) + uip(uHalf) - 2 * uHalf! \\
       13 & \lstinline!dt / dx**2 * (im(uHalf) + ip(uHalf) - 2 * uHalf)! \\
       14 & \lstinline!u0 + dt / dx**2 * (im(uHalf) + ip(uHalf) - 2 * uHalf)! \\
        \hline
    \end{tabular}
    \caption{Atomic decomposition of the midpoint scheme for the heat equation
             represented as a computational graph.}
    \label{f:decomp_heat}
\end{figure}

To algorithmically find a decomposition for a given stencil update formula,
we view it as a directed acyclic graph $(V,E)$.
Vertices in the graph, denoted by integers 0,1,$\ldots$,
represent intermediate values in the update formula.
We interchangeably use Value $i$ and Vertex $i\in V$ in this paper.
An edge $(i,j)$ exists if $j$ directly depends on $i$, i.e., if
Value $i$ is directly used in the computation of Value $j$.
Source vertices with no incoming edges are the inputs of the stencil update
formula. Sink vertices with no outgoing edges are the outputs.
Figure \ref{f:decomp_heat} shows an example of such a graph for Update formula
(\ref{update01}).
Note that each node of the computational
graph represents a symbolic value located at all the grid points, as opposed
to a value at a particular grid point.  We call this directed acyclic graph the
computational graph of the update formula.  Similar computational graphs
have been used in solving other combinatorial                   
problems in the realm of optimizing complex stencil
computation~\cite{gysi2015modesto,christen2011patus,datta2008stencil}.

Special edges represent value dependency at neighboring grid points.
If Value $j$ at a grid
point depends on Value $i$ at a neighboring grid point, then
$(i,j)$ is called a `{\bf swept}' edge.  Operations that
create swept edges include \lstinline!im(u)! and \lstinline!ip(u)!.
The set of swept edges are denoted by $E_S \subset E$.  Swept edges are
visualized by triple lines in Figure \ref{f:decomp_heat}.

The computational graph helps to visualize not only the update formula
but also its decomposition into atomic stages.  In the decomposition,
the computational graph is divided into a sequence of subgraphs,
each representing an atomic stage.  An intermediate value
can either live within a single atomic stage, or be created in one stage
and passed to subsequent stages for further use; therefore,
a vertex in the computational graph can belong to either one subgraph
or several successive subgraphs.
We do not allow a stage to repeat a computation that has been performed
in a previous stage; therefore, each vertex is created in one and only one
stage, a stage that not only solely owns all the incoming edges of that vertex,
but also contains all the vertices from which these edges originate.

This decomposition differs from many classical parallel computing methodologies,
which often starts by identifying closely related computational work and assigning them
to be executed by a single processor.  In contrast, each decomposed stage here corresponds
to a subset of computation that needs to be performed on each grid point before the
next stage can occur at the same grid point.
Neither work or data is assigned to any particular processor during this decomposition.
Each decomposed stage is atomic, and thus can be relatively easy to parallelize.

The {\bf source} of each atomic stage is the set of vertices
with no incoming edges in the subgraph; these are the inputs of the
atomic stage. The {\bf sink} of each atomic stage is the set of
vertices with no outgoing edges in the subgraph; they are the outputs.
Because the outputs of each state serve as inputs of the next, the sink of
one subgraph must be identical to the source of the next subgraph.
Also, the source of the first subgraph and the sink of the
last subgraph should match the source and sink of the entire computational
graph.

Recall that a stage is atomic only if its outputs depend on its
inputs at no further than the immediately
neighboring grid points.  This property can be graph-theoretically enforced,
by allowing {\bf at most one swept edge in any path} within a subgraph.
If every path from a source value to a sink value of the subgraph
contains no more than one swept edge, the corresponding output value can only
depend on the corresponding input value at the immediately neighboring
grid points.

Now we can formulate a graph-theoretical equivalence to the problem of
decomposing a stencil update formula into atomic stages.  This problem is
to decompose a directed acyclic graph $(V,E)$, with swept edges
$E_S\subset E$, into a sequence of subgraphs $(V_1,E_1), \ldots,
(V_{k}, E_{k})$, such that the three axioms hold:
\begin{enumerate}
    \item Each edge in $E$ belongs to one and only one subgraph.  That implies
          no redundant computation is performed.  Also, for
          each vertex $i\in V$, all incoming edges belong to one subgraph,
          which corresponds to the stage in which Value $i$ is computed.
    \item The source of each subgraph, other than the first one,
          must be contained in the previous subgraph.
          This disallows communication between non-subsequent stages.
          Also, the source of the first subgraph $(V_1,E_1)$ must
          match the source of $(V,E)$; the sink of $(V,E)$
          must be contained in $V_K$.
    \item There is no path within any subgraph $(V_k,E_k)$ that contains
          two edges in $E_S$.  This ensures that each stage is atomic.
\end{enumerate}

Figure \ref{f:decomp_heat} shows an example of an atomic decomposition.
The source of the stencil update formula is Value 0; the sink is Value 14.
Blue and green colors represent edges in the first and second subgraphs.
The source of blue subgraph includes only Vertex 0.
Vertices 1 and 3 belong exclusively to the blue subgraph; they are created in
the blue subgraph and are used to compute values only in the same
subgraph.  In other words, all their incoming and outgoing edges are blue.
Vertices 2 and 5-13 belong exclusively to the green subgraph;
all their incoming and outgoing edges are green.
Vertices 0 and 5 are shared by both subgraphs.
Both are created in the blue subgraph (inputs to the entire update formula
are defined to be created in the first subgraph); both are used in
blue and green subgraphs.  These two vertices are the sink of the blue
subgraph, and source of the green subgraph.
Neither the blue nor green subgraph contains directed path
that goes through more than one {\it swept} edges, visualized by the
triple-lines.

\section{Algebraic representation of an atomic decomposition}
\label{sec_alge}

Using the graph theoretical representation, we can formulate
a set of algebraic constraints, the satisfaction of which leads to
an atomic decomposition of an update formula.  We can then combine these
constraints with an optimization problem to
minimize the values passed between the decomposed stages.
If satisfaction of the constraints is not only sufficient but also necessary
for a valid decomposition, then by solving the constrained optimization
problem, we are guaranteed to obtain the best possible decomposition.

The primary challenge of constructing these constraints is the third criterion
in the last section, which forbids any path within a subgraph that contains
two swept edges.  Naively enforcing this criterion requires enumerating all
paths, which can be combinatorially many.  This section shows that this
criterion can be applied more efficiently, along with other criteria,
by introducing three integers for each vertex, and prescribing linear
equalities and inequalities on these integers.

To describe the set of constraints,
we use integers $1,2,\ldots,K$ to denote the $K$ decomposed stages.
We then introduce the following three integers associated with each vertex
in the computational graph:
\begin{enumerate}
    \item The {\it creating stage} $c_i$,
    \item The {\it discarding stage} $d_i$,
    \item The {\it effective stage} $e_i$,
\end{enumerate}
In this section, we first introduce $c_i, d_i$ and their governing
constraints.  Before we proceed to introduce $e_i$,
we prove the Atomic Stage Lemma, which useful for explaining
$e_i$.  Finally, we present the Quarkflow Theorem,
which lists all the constraints between $c_i, d_i,$ and $e_i$,
and shows that they are equivalent to the three criteria listed
in the previous section.

\subsection{$c_i, d_i$, and their governing constraints}
The first integer, the {\it creating stage} $c_i$,
indicates in which stage Value $i$ is created.
It is the first subgraph to which vertex $i$ belongs.
In other words, $$c_i := \min \{ k: i\in V_k\}.$$
The second integer, the {\it discarding stage} $d_i$,
indicates in which stage Value $i$ is last used.
It is the last subgraph to which vertex $i$ belongs.
In other words, $$d_i := \max \{ k: i\in V_k\}.$$

Because
a value must first be created before it is discarded, $d_i$ is always
greater or equal to $c_i$.  This leads to our first constraint,
\begin{equation} \label{cond1}
    c_i - d_i \le 0\;.
\end{equation}
Specifically, if $c_i=d_i$, then
Vertex $i$ belongs exclusively to Stage $c_i$; its corresponding vertex
has all incoming and outgoing edges in this same stage.
If $c_i<d_i$, then Vertex $i$ belongs to Stages $c_i, c_i+1, \ldots, d_i$.

If Vertex $i$ belongs to the source of the update formula, then
\begin{equation}
    c_i = 1\;,
\end{equation}
i.e., it is created at the first stage.
If Vertex $i$ belongs to the sink of the update formula, then
\begin{equation}
    d_i = K\;;
\end{equation}
this is to indicate that an output is used after the last stage.
These are equality constraints associated with the source and sink of
the update formula.

The next set of constraints is based on the following truth.
For a valid decomposition satisfying the properties in Section 3,
the creating stage satisfies the following property:

{\bf Edge Lemma:}  
For a valid decomposition,
if $(i,j)\in E_k$, then $c_j:=\min \{ k': j\in V_{k'}\}=k$.

\begin{proof} Because $(i,j)\in E_k$, $j \in V_k$.
Then by its definition, $c_j \le k$.
It is then sufficient to prove that vertex $j$ cannot be in any $V_{k'}, k'<k$.
We prove this by contradiction.  Because $j$ has one incoming edge in $E_k$,
all its incoming edges must be exclusively in $E_k$.  If $j\in V_{k'}, k'<k$,
then $j$ has no incoming edge in $(V_{k'}, E_{k'})$, and must be a source
of the subgraph.  If $k'=1$, then by Criterion 2 of Section 3, $j$ must
be in the source of $(V,E)$, which cannot be true because $j$ has an
incoming edge $(i,j)$.  If $k'>1$, then $j$ being in the source of 
$(V_{k'}, E_{k'})$ implies, by Criterion 2 again, that $j$ must be in
$(V_{k'-1}, E_{k'-1})$, which leads to contradiction by induction.
\end{proof}

{\bf Edge Corollary:}  
For a valid decomposition,
$\forall (i,j)\in E$, then $c_i:=\min \{ k': i\in V_{k'}\}\le c_j$.

\begin{proof} Each $(i,j)$ belongs to one and only one subgraph.
Denote $(i,j)\in E_k$.  Then $i\in V_k$, and by its definition $c_i\le k$.
By the Edge Lemma, $c_j=k$.  Thus $c_i\le c_j$.
\end{proof}

The lemma and corollary indicate that if $(i,j)\in E_k$, then
$c_i\le c_j=k$.  In addition, by the definition of a subgraph,
both $i$ and $j$ must belong to $V_k$.  Thus the definition of $d_i$
implies that $d_i\ge k$.  This leads to our second set of inequality
constraints,
\begin{equation} \label{cond4}
    c_i \le c_j \le d_i
\end{equation}

\subsection{The Atomic Stage Lemma}

The next set of constraints is derived from Criterion 3 of Section 3,
which enforces that the decomposed stages are atomic.  Recall that a stage
is atomic if every directed path in its subgraph contains at most one
swept edge.  Naively enforcing this criterion requires enumerating
over all paths in a subgraph, leading to exponentially many constraints.
To avoid the combinatorial explosion of constraints, we prove an equivalent
definition of an atomic stage that requires fewer constraints to enforce.

{\bf Atomic Stage Lemma:} 
A directed acyclic graph $(V_k, E_k)$ with swept edges $E_{S,k}\subset E_k$
is atomic if and only if there exists an $s_{i;k} \in \{0,1\}$ for each
$i\in V_k$,
such that
\begin{enumerate}
    \item $\forall (i,j)\in E_k$, $s_{i;k} \le s_{j;k}$,
    \item $\forall (i,j)\in E_{S;k}$, $s_{i;k} < s_{j;k}$.
\end{enumerate}

This lemma enables us to efficiently
enforce the atomicity of each stage by reducing the necessary
constraints from as many as the number of paths, which can scale exponentially
to the size of a graph, to the number of edges, which scales linearly.
This lemma makes it possible to enforce Criterion 3 in Section 3
at a computational cost that does not scale exponentially with respect to the
complexity of the stencil update formula.

\subsubsection*{Proof of the Atomic Stage Lemma}

Let us denote the directed acyclic graph $(V,E)$ with swept edge set $E_S$ as
$(V,E,E_S)$.  To prove the theorem, we decompose it into two propositions,
such that it is sufficient to prove both propositions.

{\bf Proposition 1:}
If $s_i \in \{0,1\}$ exists for each $i\in V$ and both conditions in the
theorem are satisfied, then $(V, E, E_S)$ is atomic.

{\bf Proposition 2:}
If $(V, E, E_S)$ is atomic, then there exists an $s_i \in \{0,1\}$
for each $i\in V$, such that both conditions in the
theorem are satisfied.

\begin{proof}
Proof of Proposition 1 by contradiction: If $(V,E,E_S)$ is not atomic,
then there exists a path $(i_0, \ldots, i_n)$ that contains two swept edges.
Because we can truncate both ends of such a path such that the first and last
edges are swept, we can assume $(i_0,i_1)\in E_S$ and
$(i_{n-1},i_n)\in E_S$ without loss of generality.
From Condition 2 of the theorem,
$(i_0,i_1)\in E_S\Rightarrow s_{i_0} < s_{i_1}$.  So $s_{i_0} = 0$ and
$s_{i_1} = 1$ because both must either be 0 or 1.  Similarly,
$(i_{n-1},i_n)\in E_S\Rightarrow s_{i_{n-1}} < s_{i_n}$.  So $s_{i_{n-1}} = 0$
and $s_{i_n} = 1$.  From Condition 1 of the theorem, however,
$(i_k,i_{k+1})\in E\Rightarrow s_{i_k} \le s_{i_{k+1}}$, thus the series
$s_{i_k}$ must be monotonically non-decreasing, contradicting the previous
conclusion that $s_{i_1} = 1$ and $s_{i_{n-1}} = 0$. 
\end{proof}

To prove of Proposition 2. We first construct a non-negative integer for each
vertex with the following recursive formula
\begin{equation} \label{proof_construct_s_j}
    s_j = \begin{cases}
        0, & \mbox{Vertex j has no incoming edge} \\
        \max_{i : (i,j)\in E'} s_i + I_{(i,j)\in E_S} , & \mbox{otherwise}
    \end{cases}
\end{equation}
where $I_{(i,j)\in E_S} = 1$ if $(i,j)\in E_S$ and 0 otherwise.
By construction, these integers satisfy both conditions in the theorem.
We then show that if the graph is atomic, then all these integers
satisfy the additional condition that $s_i \in \{0,1\}$.
This additional condition can be proved with the aid of the following lemma:

{\bf Lemma:} For every vertex $j$, there exists a path, consisting of zero
or more edges, that ends at $j$ and contains $s_j$ swept edges, where
$s_j$ is defined by Equation (\ref{proof_construct_s_j}).

    \begin{proof}
Proof of Lemma by induction:
We prove by induction with respect to the length of the longest path that
ends at Vertex $j$.  If the length of the longest incoming path is 0,
it means Vertex $j$ has no incoming edge; by definition, $s_j=0$.
In this case, an empty path suffices
as one that ends at $j$ and contains $0$ swept edges.  Therefore, the lemma
is true if the longest incoming path to $j$ is of length 0.

If the lemma is true for any vertex whose longest incoming path is of length
less than $n>0$, we prove it for any vertex $j$ whose longest incoming path
is of length $n$.  Because $n>0$, Vertex $j$ has an incoming edge.
By Equation (\ref{proof_construct_s_j}), there exists an edge $(i,j)$ such that
$s_j = s_i + I_{(i,j)\in E_S}$.  Note that Vertex $i$ has no incoming path
of length $n$; if it does, then appending edge $(i,j)$ to that path would
result in a path to Vertex $j$ of length $n+1$, violating the assumption
about the longest incoming path to $j$.  Because the longest
incoming path to Vertex $i$ is less than $n$, by the induction assumption,
the lemma holds for Vertex $i$:
there exists a path that ends at Vertex $i$ that contains $s_i$ swept edges.
At the end of this path, we append edge $(i,j)$ to form a path to Vertex
$j$.  If $(i,j)\notin E_S$, the path contains the same number of swept edges
as the path to $i$.  $(i,j)\notin E_S$ also sets $s_j=s_i$
by Equation (\ref{proof_construct_s_j}).
Thus this path suffices as one that ends at Vertex $j$ and contains $s_j$
swept edges.  The other case is $(i,j)\in E_S$,
and the resulting path to $j$ has one more swept edge than the path to $i$.
In this case, $s_j=s_i+1$ according to Equation (\ref{proof_construct_s_j}),
and this path still suffices as one that ends at Vertex $j$ and contains $s_j$
swept edges.  This shows that the lemma is true for any Vertex $j$ whose
longest incoming path is of length $n$, thereby completing the induction.
    \end{proof}

    \begin{proof}
Proof of Proposition 2 by contradiction:
If the proposition is not true for an atomic stage,
then $s_j$ constructed via Equation 
(\ref{proof_construct_s_j}) must not satisfy $s_j\in \{0,1\}$ for some $j$.
Because $s_j$ are non-negative integers by construction, 
$s_j\notin \{0,1\}$ means $s_j>1$.  By the lemma, there exists a path that
ends at Vertex $j$ that contains more than 1 swept edges. So the graph
cannot be an atomic stage, violating the assumption. 
    \end{proof}

\subsection{The Quarkflow Theorem}
The set of integers $s_{i;k}$ in the lemma, however, is specific to each
subgraph $k$.
They can only be defined on top of a valid decomposition. 
To specify a-priori conditions that can serve as foundation to build a
valid decomposition, we construct $e_i$, the effective stage,
for each $i\in V$.  The construction of $e_i$ is such that for
each $0\le k<K$,
\begin{equation} \label{s_i_e_i}
    s_{i;k} := \begin{cases}
        e_i - c_i \;, & c_i = k \\
        0\;, & c_i < k
    \end{cases}
\end{equation}
satisfies both conditions in the Atomic Stage Lemma for subgraph $k$.
The introduction of $e_i$ and its associated constraints,
together with $c_i$ and $d_i$ and their associated constraints derived
previously, leads us to the main theorem of this section.
We name it the ``Quarkflow``
theorem since the triplet of integers $(c_i,d_i,e_i)$ representing each
vertex is analogous to quarks forming each subatomic particle.

{\bf Quarkflow Theorem:} 
A graph $(V,E)$ is decomposed into $K$ subgraphs,\\
$(V_0,E_0)$, $\ldots$, $(V_{K-1},E_{K-1})$.
The decomposition satisfies the three criteria stated
in Section 3, if and only if there exists a triplet of integers
$(c_i, d_i, e_i)$ for vertex $i\in V$, such that the following constraints
are satisfied:
\begin{enumerate}
    \item $c_i \le d_i$ for all $i\in V$
    \item $c_i = 1$ for all $i$ in the source of $V$
    \item $d_i = K$ for all $i$ in the sink of $V$
    \item $c_i \le c_j \le d_i$ for all $(i,j)\in E$
    \item $c_i \le e_i \le c_i + 1$ for all $i\in V$
    \item $e_i \le e_j$ for all $(i,j)\in E$
    \item $e_i + 1 \le e_j$ for all $(i,j)\in E_S$
    \item $c_i + 1 \le e_i$ for all $i\in V$ where $\exists (j,i)\in E_S$
\end{enumerate}
and that $\forall k=1,\ldots,K,
V_k=\{i : c_i\le k\le d_i\}$, $E_k=\{(i,j)\in E : c_j=k\}$.

This theorem achieves the goal of this section. To find a decomposition,
we only need to find the integer triplets, $(c_i,d_i,e_i),\forall i\in V$,
together with the number of stages $K$, satisfying all these constraints.
If we can use the same set of integers
to describe how good the decomposition is, we can solve an integer program
to find the best possible decomposition.

\subsubsection*{Proof of the Quarkflow Theorem}

We split the proof into two subsections.  The first shows that 
for a decomposition that satisfies the three criteria in Section 3,
there exists $(c_i,d_i,e_i),\forall i\in V$ satisfying
all constraints in the theorem, and the decomposition matches $(V_k,E_k)$
constructed in the theorem.  The second subsection shows that
if the triplets satisfy all the constraints in the theorem, then 
$(V_k,E_k)$ is a valid decomposition satisfying all three criteria
in Section 3.

\subsubsection{Criteria in Section 3 $\Longrightarrow$
Constraints in the Quarkflow Theorem}
For a decomposition that satisfies the three criteria in Section 3,
we can construct $(c_i,d_i,e_i)$ in the following way:
\begin{enumerate}
    \item $c_i = \min \{ k | i\in V_k \}$
    \item $d_i = \max \{ k | i\in V_k \}$
    \item $e_i = c_i + s_{i;c_i}$, where \{$s_{i;c_i}\}_{i\in V_{c_i}}$
        satisfy the two constraints in the Atomic
        Stage Lemma for subgraph $c_i$.
\end{enumerate}
Constraints 1, 2, and 3 are satisfied by construction.

For Constraint 4, consider each $(i,j)\in E$ belonging to
subgraph $E_k$.  Therefore, $i\in V_k$ and $j\in V_k$.
Thus, $c_i\le k$ and $d_j\ge k$.  By the Edge Lemma, $c_j=k$.
Constraint 4 therefore holds.

Constraint 5 is satisfied by construction of $e_i$ and by the Atomic Stage
Theorem.

Constraint 6 and 7 can be proved in two cases.
The first case is when $c_i=c_j$.  Denote both $c_i$ and $c_j$ as $k$.
From the Atomic Stage Lemma,
$s_{i;k} \le s_{j;k}$, and when $(i,j)\in E_S$, $s_{i;k}+1 \le s_{j;k}$.
Plug these inequalities, as well as $c_i=c_j=k$, into the definition
of $e_i$ and $e_j$, we get
$e_i \le e_j$ and when $(i,j)\in E_S$, $s_i+1 \le s_j$.

The second case is when $c_i<c_j$ (the Edge Corollary ensures that $c_i$
cannot be greater than $c_j$).  Because both $c_i$ and $c_j$
are integers, $c_i+1\le c_j$.  Also because both $s_{i;c_i}$ and $s_{j;c_j}$
are either 0 or 1,
$e_i := c_i+s_{i;c_i}\le c_i \le c_j \le c_j+s_{j;c_j}=:e_j$.
In particular, when $(i,j)\in E_S$, the Atomic Stage Lemma ensures that
$s_{i;c_j} < s_{j;c_j}$; since both are either 0 or 1, the strict inequality
leads to $s_{j;c_j}=1$.  Thus
$e_i := c_i+s_{i;c_i}\le c_i+1 \le c_j = c_j+s_{j;c_j}-1=:e_j-1$, and thus
$e_i+1\le e_j$.

To prove Constraint 8 for an edge $(i,j)\in E_S$, we use Criterion 1 of
Section 3, which specifies that $(i,j)$ is in a unique subgraph, denoted
as $E_k$.  By the Edge Lemma, $k=c_j$.  Because $(i,j)\in E_S\cap E_k$,
the Atomic Stage Lemma guarantees that $s_{i;k} < s_{j;k}$.  Because both
$s_{i;k}$ and $s_{j;k}$ are 0 or 1, this strict inequality ensures
that $s_{j;k}=1$.  Therefore, $e_j:= c_j + s_{j;k} \ge c_j+1$, which is
Constraint 8.

Finally, we show that $V_k=\{i : c_i\le k\le d_i\}$ and
$E_k=\{(i,j)\in E : c_j=k\}$.  For $V_k=\{i : c_i\le k\le d_i\}$ to be
true, it is sufficient to show that $i\in V_k\Leftrightarrow c_i\le k\le d_i$.
If $c_i>k$ or $d_i<k$, then $i\notin V_k$
by the construction of $c_i$ and $d_i$.  It is then sufficient to
show that $i\in V_k$ for all $k$ satisfying $c_i\le k\le d_i$.  We prove
this by contradiction.  If the previous statement is not true, then
because $i\in V_{c_i}$ and $i\in V_{d_i}$ by the construction of $c_i$
and $d_i$, there must exist $c_i < k < d_i$ such that
$i\in V_{k+1}$ but $i\notin V_k$.  By the edge lemma, all incoming edges
of $i$ must be in subgraph $c_i$, so $i$ must be in the source of subgraph
$k+1$.  By Criterion 2 of Section 3, $i\in V_k$ which contradicts with
$i\notin V_k$.  Therefore $i\in V_k\Leftrightarrow c_i\le k\le d_i$ and
$V_k=\{i : c_i\le k\le d_i\}$.

For $E_k=\{(i,j)\in E : c_j=k\}$, we know from the edge lemma that
$(i,j)\in E_k \Rightarrow c_j=k$.  From Criterion 1, each edge can only be
in one subgraph.  Therefore, $(i,j)\in E_k \Leftarrow c_j=k$.
Thus $E_k=\{(i,j)\in E : c_j=k\}$. \qed

\subsubsection{Constraints in the Quarkflow Theorem $\Longrightarrow$
Criteria in Section 3}

Now if all 8 constraints are satisfied, we prove that 
$(V_k:=\{i : c_i\le k\le d_i\},E_k:=\{(i,j)\in E : c_j=k\})$
is a graph decomposition satisfying all three criteria specified in Section 3.

First, we need to show that each $(V_k,E_k)$ is a subgraph.
If $(i,j)\in E_k$, then by the definition of $E_k$ and
Constraint 1 of the Quarkflow Theorem,
$k=c_j\le d_j$.  This ensures that $j\in V_k$ by the definition of $V_k$.
On the other hand, Constraint 4 of the Quarkflow
Theorem ensures that $c_i\le c_j\le d_i$.  Because $c_j=k$,
$c_i\le k\le d_i$.  This ensures that $i\in V_k$ by the definition of $V_k$.
Because $i\in V_k$ and $j\in V_k$ for every $(i,j)\in E_k$,
$(V_k,E_k)$ is a valid subgraph.

Now we proceed to prove that the subgraphs satisfy the three criteria
of Section 3.

{\bf Criterion 1:}
Each edge $(i,j)$ belongs to Subgraph $c_j$ and no other subgraph by the
definition of $E_k$.  Also by this definition, for each vertex $j\in V$,
all incoming edges belongs to Subgraph $c_j$.  Criterion 1 therefore holds.

{\bf Criterion 2:}
This criterion has three statements.
We first prove the first statement:
if vertex $j$ is in the source of Subgraph $(V_k,E_k), k>1$,
then in $j\in V_{k-1}$.

By the definition of $V_k$, $c_j\le k\le d_j$.
We now prove that $c_j<k$.  This is because if $c_j=k>1$,
then by Constraint 2 of the Quarkflow Theorem,
$j$ is not in the source of $V$, i.e.,
$\exists (i,j)\in E$.  $c_j=k$ also implies,
by the definition of $E_k$, that such $(i,j)\in E_k$.  This contradicts
with the assumption that $j$ is in the source of $(V_k,E_k)$.

Now we know that $c_j< k\le d_j$, which implies that
$c_j\le k-1< d_j$.  Therefore, $j\in V_{k-1}$ by the
definition of $j\in V_{k-1}$.

We then prove the second statement of the criterion:
$j$ is in the source of $(V_1,E_1) \Longleftrightarrow$ it is in the source of $(V,E)$.

$j$ being in the source of $(V_1,E_1)$ implies that $j\in V_1$, and by
the definition of $V_1$, $c_j = 1$.  By the definition of $E_1$,
any $(i,j)\in E$ must be in $E_1$.  Such incoming edge
$(i,j)$ cannot exist if $j$ is in the source of $(V_1,E_1)$.  $j$,
therefore, has no incoming edge and is in the source of $(V,E)$.

On the other hand, by Constraint 2 of the Quarkflow Theorem,
$j$ being in the source of $(V,E)$ implies that $c_j=1$.
By Constraint 1 of the Quarkflow Theorem, $c_j=1\le d_j$,
which means $j\in V_1$ by the definition of $V_1$.
$j$ has no incoming edge at all, and
therefore is in the source of $(V_1,E_1)$.

We finally prove the third statement of the criterion.
if $j$ is in the sink of $(V,E)$, then by Constraint 3 of the Quarkflow
Theorem, $d_j=K$.  By Constraint 1 of the Quarkflow Theorem,
$c_j\le K=d_j$.  Therefore, by the definition of $V_k$, $j\in V_k$.
This completes the proof of Criterion 2 of Section 3.

{\bf Criterion 3:} To prove Criterion 3, we use
Equation (\ref{s_i_e_i}) to construct
a set of $s_{i,k}$ for each subgraph $k$ and $i\in V_k$.
We now show that such $s_{i,k}$ satisfies both conditions
in the Atomic Stage Lemma. Thus Criterion 3 holds by the Lemma.

To show this, consider any $(i,j)\in E_k$.  By the construction of
$E_k$, we have $c_j = k$.
By Constraint 4 of the Quarkflow Theorem, $c_i \le c_j=k$.
The rest of the proof splits into two cases, $c_i = c_j=k$ and $c_i < c_j=k$.
\begin{itemize}
\item If $c_i=c_j=k$, then by Equation (\ref{s_i_e_i}),
$$s_{j;k} := e_j-c_j=e_j-k, \quad
  s_{i;k} := e_i-c_i=e_i-k.$$
Constraint 6 and 7 of the Quarkflow Theorem states that
$$e_j - e_i \ge
\begin{cases}
    1, & (i,j)\in E_S \\
    0, & (i,j)\notin E_S
\end{cases}.$$
Therefore,
$$ s_{j;k} - s_{i;k} = (e_j-k)-(e_i-k) = e_j-e_i\ge
\begin{cases}
    1, & (i,j)\in E_S \\
    0, & (i,j)\notin E_S
\end{cases},$$
which is a sufficient condition for the Atomic Stage Lemma to hold.
\item If $c_i< c_j=k$, then by Equation (\ref{s_i_e_i}),
$$s_{j;k} := e_j-c_j, \quad
  s_{i;k} := 0.$$
So, $s_{j;k} - s_{i;k} = e_j-c_j$.
Using Constraints 5 and 8 of the Quarkflow Theorem, we therefore have
$$ s_{j;k} - s_{i;k} \ge
\begin{cases}
    1, & (i,j)\in E_S \\
    0, & (i,j)\notin E_S
\end{cases}.$$
\end{itemize}
which is sufficient for the Atomic Stage Lemma.
In both cases, we now have proved Criterion 3. \qed

\section{Decomposition as the dual of a network flow problem}

In Section \ref{sec_graph}, we formulated three criteria that define a valid
decomposition of a stencil update formula into atomic stages.
In Section \ref{sec_alge}, we proved that any such valid decomposition
could be represented by a set of integers, $K$
and $(c_i,d_i,e_i),\forall i\in V$ that satisfy a set of linear constraints
among them.  The reverse is also true: any set of integers satisfying
these constraints represents a valid decomposition.  In this section, we
complete these linear constraints with a linear cost function to form
an integer program.

The cost function should model how expensive it is to execute the
decomposed formula.  This cost, in a massively parallel computation, depends
on the interconnect, the domain decomposition, and the solver.
In this paper, we use a very simplistic cost function, composed of
two factors:
the number of stages and the amount of coupling between the stages.
These factors correspond to the communication cost
of executing the decomposed formula in the swept rule of parallel
computing~\cite{alhubail2016swept}.  The number of stages is proportional to how often
data needs to be communicated, and thereby models the communication
time due to network latency.  The amount of coupling between stages can
correspond to how much data is communicated, and thereby models
how much time is spent due to network bandwidth.
Both factors can be represented algebraically using the same
integers $K$ and $(c_i,d_i,e_i)$ we use in Section 3.

The first factor, the number of atomic stages, is simply $K$.  The second
factor, the amount of coupling, can be modeled as the amount of
values shared between stages.  If a value belongs to only one stage,
i.e., $c_i=d_i$, it does not contribute to the amount of coupling.
If a value belongs to multiple stages, it couples these stages $c_i,\ldots,d_i$.
It then adds an amount $w_i (d_i-c_i)$ to the total coupling.
$w_i,i\in V$ here are the weights of the vertices in the graph.  If all vertices
represent values of the same size, we can set $w_i\equiv 1$.  If vertices
can represent values of different size, i.e., vectors, matrices, and tensors,
then $w_i$ can be adjusted to the number of bytes required to store the $i$th
intermediate value at each grid point.  The two factors can be blended using
a positive weight factor $W_K$, the only tuning parameter of the cost function.
Curiously, the choice of $W_K$
does not appear to affect the outcome of the optimization in the test cases
we have so far examined.  The blended cost function combines with the
constraints in the Quarkflow Theorem of Section 4 to form the following integer
programming problem,
\begin{equation} \label{ip} \begin{split}
    \min_{K,c_i,d_i,e_i,i\in V} \quad & W_K K + \sum_i w_i (d_i - c_i) \\
    s.t.  \quad & K,c_i,d_i,e_i\in \mathcal{Z},\quad \mbox{and}\\
    & \mbox{all 8 conditions in the Quarkflow Theorem are satisfied.}
\end{split} \end{equation}
Here $\mathcal{Z}$ represent integers.

This integer program problem can be solved by solving the equivalent linear
program, ignoring the requirement that the solutions must be integers,
by a simplex method.  This is a characteristic shared by all integer
programming problems in which the constraint matrix is totally
unimodular~\cite{schrijver1998theory} and the constraint constants are integers.
This problem (\ref{ip}) has a unimodular constraint matrix because all its
constraints involve the difference between pairs of variables~\cite{garfinkel1972integer}.
Therefore, each row of the constraint matrix has exactly two entries; one
equals to 1, the other equals to -1.  The transpose of this matrices
has the same property as the constraint matrix of a flow network.
Because the constraint matrix of flow networks are totally
unimodular~\cite{ahuja1993network}, our constraint matrix is totally unimodular.

Not only is the constraint matrix transpose to that of a flow network,
the linear program is the symmetric dual of
a minimum-cost flow problem on a flow network.  In this flow network,
each variable in Problem (\ref{ip}), including $c_i, d_i, e_i$ and $K$,
is a vertex.  An additional vertex, denoted as 0,
is introduced for the single-variable constraints
$c_i = 1,\forall i \in V_S$.  The dual of Problem (\ref{ip}) as
a linear program is
\begin{equation} \label{dual_lp} \begin{split}
    \mbox{Choose} \; & \\
    \forall i\in V_S, & \quad x_{c_i,0}\in R, \\
    \forall i\in V_T, & \quad x_{d_i,K}\in R, \\
    \forall i\in V,   & \quad x_{c_i,d_i}\ge 0, \\
                      & \quad x_{c_i,e_i}\ge 0, \\
                      & \quad x_{e_i,c_i}\ge 0, \\
    \forall (i,j)\in E, & \quad x_{c_i,c_j} \ge 0, \\
                        & \quad x_{c_j,d_i} \ge 0, \\
                        & \quad x_{e_i,e_j} \ge 0 \\
    \mbox{in order to} \; & \\
    \min \quad & \sum_{i\in V_S} x_{c_i,0} + \sum_{i\in V} (x_{e_i,c_i} -
                 I_{\exists (j,i)\in E_S} x_{c_i,e_i})
               - \sum_{(i,j)\in E_S} x_{e_i, e_j} \\
    s.t. \quad
    & \sum_{i\in V_T} x_{d_i,K} = W_K, \quad
    \mbox{and } \forall i\in V, \\
    & I_{i\in V_S} x_{c_i,0} +
      x_{c_i, d_i} + x_{c_i, e_i} - x_{e_i, c_i} \\
    &+ \sum_{j: (j,i)\in E} (x_{c_i, d_j} - x_{c_j, c_i}) +
      \sum_{j: (i,j)\in E} x_{c_i, c_j} = w_i \\
    & I_{i\in V_T} x_{d_i, K} - x_{c_i, d_i} -
      \sum_{j: (i,j)\in E} x_{c_j, d_i} = -w_i \\
    & -x_{c_i, e_i} + x_{e_i, c_i} + \sum_{j:(i,j)\in E} x_{e_i,e_j}
      - \sum_{j:(j,i)\in E} x_{e_j,e_i} = 0 \\
\end{split} \end{equation}
The constraints are flow balance equations on all vertices in the flow network,
except for vertex 0, whose balance equation is a redundant linear combination
of the constraints listed above:
\[ \sum_{i\in V_S} x_{c_i,0} = -W_K \]
Adding this redundant constraint to the above problem, it becomes clear
that this is a minimum cost flow network problem~\cite{ahuja1988network},
one that can leverage specialized, fast solvers, such as the out-of-kilter
algorithm~\cite{fulkerson1961out, ahuja1993network}
and relaxation algorithms~\cite{bertsekas1988relaxation}.

\section{Implementation, Tests and Results}

To demonstrate the algorithm developed in this paper,
the minimum cost flow network problem (\ref{dual_lp}) is
solved using the out-of-kilter algorithm~\cite{fulkerson1961out},
implemented in the open source package GLPK~\cite{makhorin2008glpk}.
The out-of-kilter algorithm returns the
solution of both the linear program (\ref{dual_lp}) and
its dual (\ref{ip}).  The dual solution includes $c_i$ for each vertex $i\in V$,
indicating the stage in which the corresponding value should be created.
The solution also includes $d_i$ for each vertex.  Comparing it
to $c_i$ tells us whether Value $i$ is used only in its
creating stage, or should be passed to subsequent stages for further use.

\begin{figure}[htb!] \centering
    \subfloat[]{\includegraphics[width=0.6\textwidth]{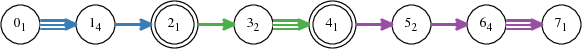}}\\
    \subfloat[]{\includegraphics[width=0.6\textwidth]{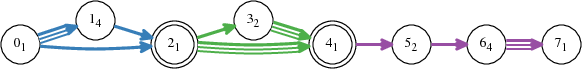}}\\
    \subfloat[]{\includegraphics[width=0.6\textwidth]{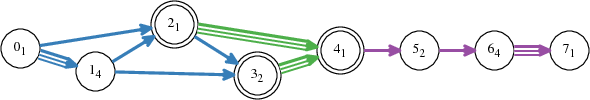}}\\
    \subfloat[]{\includegraphics[width=0.6\textwidth]{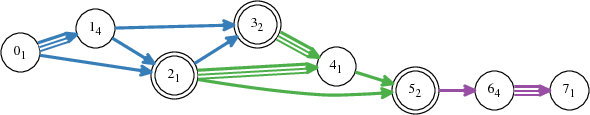}}\\
    \subfloat[]{\includegraphics[width=0.6\textwidth]{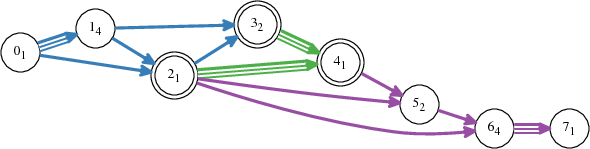}}\\
    \subfloat[]{\includegraphics[width=0.6\textwidth]{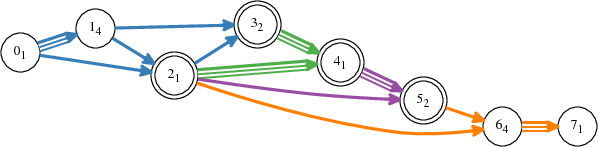}}
    \caption{Atomic decomposition of manufactured test cases.
        These graphs are not generated from real stencil update formulae
        but are constructed for testing and illustration purposes.
    The number and subscript shown in each vertex represent the index $i$
    and weight $w_i$ of the vertex.  Triple-lined edges represent ``swept''
    edges.  Different colors represent different stages.
    Double circles represent vertices shared by more than one stages.}
    \label{f:manu}
\end{figure}

The implementation can be found in the ``Quarkflow'' tool included
in the software ``pascal''\footnotemark[1]
\footnotetext[1]{https://github.com/qiqi/pascal,
branch master, commit \\df01f76a8bf1bc64237d2e346616dd67e6dd8116}.
The algorithm is first tested in a series of manufactured test cases,
as shown in Figure \ref{f:manu}.

In (a), our algorithm is asked
to decompose a graph of eight vertices; each vertex only depends on the
previous one.  Three out of the seven dependencies are swept.
Under the constraints of atomic decomposition, the
optimal decomposition should split the graph into three subgraphs,
each containing one swept edge.  It should also minimize the total
weight of the vertices shared by two subgraphs.  The outcome of 
our algorithm produces this behavior; both shared
vertices are of weight 1.

In (b), two additional edges are added, but without affecting the
decomposition.  The creating stage and discarding stage of each variable
can stay the same, thereby creating no more shared vertices.
Our algorithm still chooses a similar decomposition as before.

In (c), an edge is added from Vertex 1 to Vertex 3.  If Vertex 3 had
the same creating stage, green, as in (b), Vertex 1 must be shared
between the blue and green stages.
This sharing would increase the cost function by 4, the weight of Vertex 1.
To reduce the cost, the algorithm moves Vertex 3 into the blue stage, sharing
it between the blue and green stages.  Because Vertex 3 has a weight of 2,
sharing it would increase the cost function only by 2.

In (d), another edge is added from Vertex 2 to Vertex 5.  Keeping the creating
stage of Vertex 5 would require now sharing Vertex 2 between all three stages,
increasing the cost function of (\ref{ip}) by 1.  Instead, our algorithm
changed the creating stage of Vertex 5, sharing it between the green and
purple stages.  The resulting increase of the cost function is also 1. 
Multiple integer optimums exist in this case.
Note that any linear combination of these integer optimums also solves
the linear program.  The out-of-kilter algorithm always chooses one of
the integer optimums.

In (e), an edge is added from Vertex 2 to Vertex 6.  Now Vertex 2 must be
shared among all three stages anyway, our algorithm found it unnecessary to
share Vertex 5, and reverts to sharing the cheaper Vertex 4 between the green
and purple stages.

In (f), the edge from 4 to 5 is changed to a swept edge.  It forces an
additional stage to be created, by splitting from
the purple stage in (e), sharing Vertices 2 and 5 with the remaining purple
stage.  Sharing Vertices 2 and 3, the combined weight of which is 3, is
cheaper than the alternative splitting of
sharing Vertex 6, which alone has a weight of 4.

\begin{figure}[htb!] \centering
    \includegraphics[width=\textwidth]{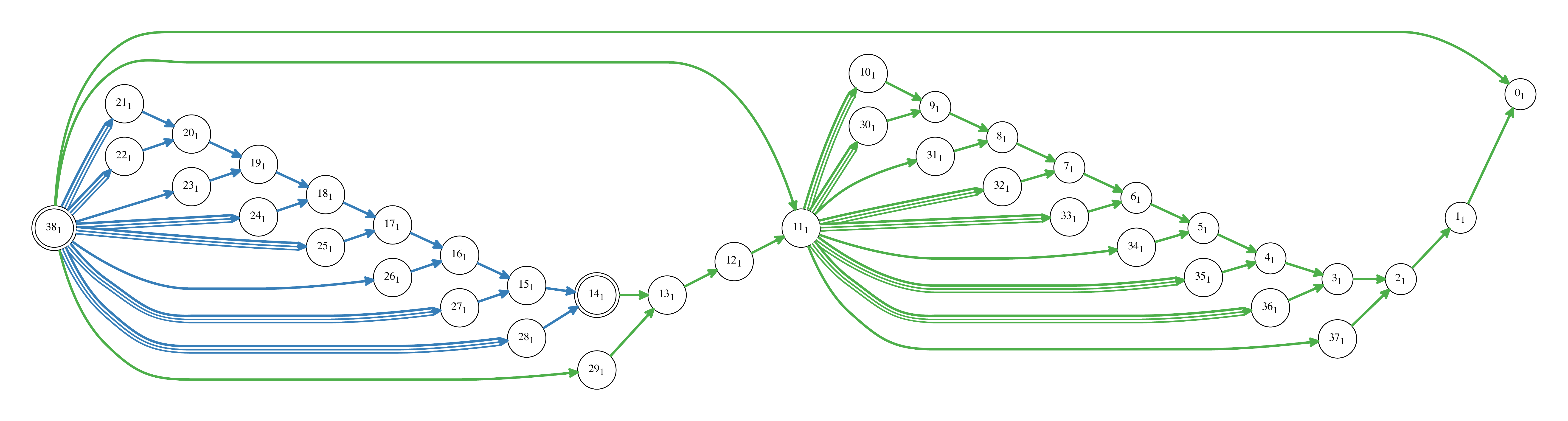}
    \caption{Atomic decomposition of a 3D heat equation, integrated using
    the midpoint rule (two-stage Runge-Kutta).
    The number and subscript shown in each vertex represent the index $i$
    and weight $w_i$ of the vertex.  Triple-lined edges represent ``swept''
    edges.  Different colors represent different stages.  Double circles
    represent vertices shared by more than one stages.}
    \label{f:heat3d}
\end{figure}

Our algorithm is then tested on a stencil update formula for solving the
3D heat equation using a two-stage Runge-Kutta scheme.
The formula is described in Python
as following.
\begin{lstlisting}
        uh = u + 0.5 * dt / dx**2 * (im(u) + ip(u) - 2 * u +
                                     jm(u) + jp(u) - 2 * u +
                                     km(u) + kp(u) - 2 * u)
        return u + dt / dx**2 * (im(uh) + ip(uh) - 2 * uh +
                                 jm(uh) + jp(uh) - 2 * uh +
                                 km(uh) + kp(uh) - 2 * uh)
\end{lstlisting}
The code above is fed into the software ``pascal''\footnotemark[1],
which automatically builds the computational graph using operator overloading,
then calls ``quarkflow'' to decompose the graph, then translates each
subgraph into C code for execution on parallel architectures.  Figure
\ref{f:heat3d} illustrates the decomposed computational graph in this
example.  The computational graph contains 39 vertices connected by 56
directed edges, among which 6 are swept.  The decomposition algorithm
completed in about 0.38 milliseconds on an Intel(R) Core(TM) i7-6650U
processor clocked at 2.20GHz.  The update formula is decomposed into
two stages, containing 22 and 34 edges, respectively.  They include
16 and 25 vertices, sharing Vertices 14 and 38.

\begin{figure}[htb!] \centering
    \includegraphics[width=\textwidth]{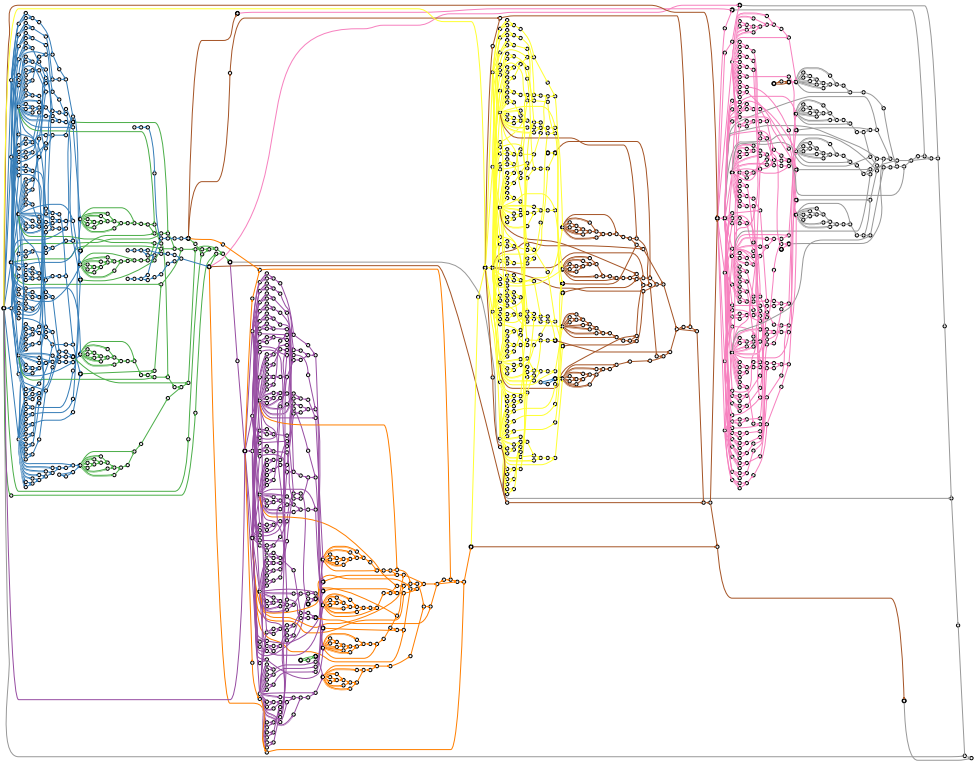}
    \caption{Automatic decomposition of a complex time step into eight atomic
    stages.  The update formula is for the Euler equation of
    gas dynamics in three spatial dimensions, using a second-order
    finite-difference spatial discretization scheme, and a fourth-order
    Runge-Kutta time integrator.
    }
    \label{f:euler3d}
\end{figure}

Finally, our algorithm is applied to an update formula for the Euler equation
of gas dynamics.  The spatial discretization, described in Appendix A,
is the conservative,
skew-symmetric finite difference scheme, with a conservative, symmetric,
negative-definite 4th order numerical dissipation. The temporal discretization
is the 4th order Runge-Kutta scheme.  Each spatial discretization requires
information from the neighbor of neighboring grid points.  Each time step
requires 4 level of spatial discretization, thereby requiring access to 8
levels of neighboring grid points.  The formula was encoded in Python below
and is automatically decomposed.  Figure \ref{f:euler3d} shows the
decomposition into 8 stages of this update formula.
The computational graph has 1424 vertices and 2106 edges, among which 432
are swept.  The decomposition algorithm
completed in about 0.25 seconds on an Intel(R) Core(TM) i7-6650U
processor clocked at 2.20GHz.  The update formula is decomposed into
eight stages, containing an average of 186.625 nodes.  An average of 6.785 nodes
in each stage is shared with subsequent stages.

\section{Conclusion}

This paper shows that we can decompose a complex stencil update
formula by solving the dual of a minimum cost flow network problem.
This is a highly nontrivial result because brute force solution of
the original problem is an integer program of combinatorial complexity.
Instead, the reformulation presented in this paper leads to a
linear program of appealing structure, amenable to some very efficient
algorithms.  As a result, even complex stencils can be decomposed
efficiently.

Efficient and automatic decomposition of complex stencil update formulas,
such as the one for the Euler equation presented in Section 6, simplifies
application of the swept decomposition rule~\cite{alhubail2016swept}.
To use the swept decomposition rule to a numerical scheme, one
only need to encode the stencil update formula of the scheme in Python.
The update formula is then quickly and automatically decomposed into a series
of atomic stages, each can be translated into low level computer code.
The resulting code can then be inserted into a software framework for executing
the scheme in the swept decomposition rule.

\section*{Acknowledgment}

The author acknowledges support from the NASA TTT program Award NNX15AU66A
under Dr. Mujeeb Malik and technical monitor Dr. Eric Nielsen, and from
the DOE Office of Science ASCR program Award DE-SC-0011089 under
Dr. Sandy Landsberg and Dr. Steve Lee.
This paper significantly benefited from discussions with Dr. David Gleich,
Dr. Christopher Maes, and Dr. Michael Saunders, as well as comments from
the anonymous reviewers.


\appendix
\section{Discretization of the 3D Euler equation}
\begin{lstlisting}
def diffx(w): return (ip(w) - im(w)) / (2 * dx)
def diffy(w): return (jp(w) - jm(w)) / (2 * dy)
def diffz(w): return (kp(w) - km(w)) / (2 * dz)

def div_dot_v_phi(v, phi):
    return diffx(v[0] * phi) + diffy(v[1] * phi) + diffz(v[2] * phi)

def v_dot_grad_phi(v, phi):
    return v[0] * diffx(phi) + v[1] * diffy(phi) + v[2] * diffz(phi)

def dissipation(r, u):
    laplace = lambda u: (ip(u) + im(u) + jp(u) + jm(u) + jp(u) + jm(u)) / 6 - u
    return laplace(DISS_COEFF * r * r * laplace(u))

def assemble_rhs(mass, momentum_x, momentum_y, momentum_z, energy):
    rhs_w = zeros(w.shape)
    rhs_w[0] = -mass
    rhs_w[1] = -momentum_x
    rhs_w[2] = -momentum_y
    rhs_w[3] = -momentum_z
    rhs_w[4] = -energy
    return rhs_w

def rhs(w):
    r, rux, ruy, ruz, p = w
    ux, uy, uz = rux / r, ruy / r, ruz / r
    ru, u = (rux, ruy, ruz), (ux, uy, uz)

    mass = div_dot_v_phi(ru, r)
    mom_x = (div_dot_v_phi(ru, rux) + r * v_dot_grad_phi(ru, ux)) / 2 + diffx(p)
    mom_y = (div_dot_v_phi(ru, ruy) + r * v_dot_grad_phi(ru, uy)) / 2 + diffy(p)
    mom_z = (div_dot_v_phi(ru, ruz) + r * v_dot_grad_phi(ru, uz)) / 2 + diffz(p)
    energy = gamma * div_dot_v_phi(u, p) - (gamma-1) * v_dot_grad_phi(u, p)

    dissipation_x = dissipation(r, ux) * c0 / dx
    dissipation_y = dissipation(r, uy) * c0 / dy
    dissipation_z = dissipation(r, uz) * c0 / dz

    rhs_r  = c0 * fan * (r  - r_fan)
    rhs_ux = c0 * fan * (ux - ux_fan) + c0 * obstacle * ux
    rhs_uy = c0 * fan * (uy - uy_fan) + c0 * obstacle * uy
    rhs_uz = c0 * fan * (uz - uz_fan) + c0 * obstacle * uz
    rhs_p  = c0 * fan * (p  - p_fan)

    return assemble_rhs((mass - rhs_r) / (2 * r) + fan * (r - r_fan),
                        (mom_x + dissipation_x - rhs_ux) / r,
                        (mom_y + dissipation_y - rhs_uy) / r,
                        (mom_z + dissipation_z - rhs_uz) / r, (energy - rhs_p))

def step(w):
    dw0 = dt * rhs(w)
    dw1 = dt * rhs(w + 0.5 * dw0)
    dw2 = dt * rhs(w + 0.5 * dw1)
    dw3 = dt * rhs(w + dw2)
    return w + (dw0 + dw3) / 6 + (dw1 + dw2) / 3
\end{lstlisting}

\bibliographystyle{siamplain}
\bibliography{master}
\end{document}